\documentclass[a4paper,
               reqno, 
               twocolumn]{article}

\usepackage{amsmath,amssymb,amsthm,mathtools}
\usepackage[utf8]{inputenc}
\usepackage{array,booktabs}
\usepackage[center,skip=4pt,font=small,labelfont=bf,width=8cm]{caption}
\usepackage{enumitem,xcolor}

\usepackage{sectsty}
\sectionfont{\large}

\providecommand{\Real}{\mathbb{R}}
\providecommand{\Integer}{\mathbb{Z}}
\providecommand{\CyI}[1]{\Integer[\zeta_{#1}]}
\providecommand{\abs}[1]{\lvert#1\rvert}

\definecolor{lgrey}{gray}{0.35}

\begin{document}

\title{\textbf{Tiling Vertices and the Spacing Distribution\\
               of their Radial Projection}}

\author{Tobias Jakobi}

\date{{\small Fakultät für Mathematik, Universität Bielefeld\\ 
              Universitätsstraße 25, D--33615 Bielefeld, Germany}}

\twocolumn[
\begin{@twocolumnfalse}
  \maketitle
  \begin{abstract} 
    \noindent The Fourier-based diffraction approach is an established method to extract order and symmetry properties
    from a given point set. We want to investigate a different method for planar sets which works in direct space
    and relies on reduction of the point set information to its angular component relative to a chosen reference
    frame. The object of interest is the distribution of the spacings of these angular components, which can for instance be
    encoded as a density function on $\Real_{+}$. In fact, this \emph{radial projection} method is not entirely new, and
    the most natural choice of a point set, the integer lattice $\Integer^2$, is already well understood.\\
    We focus on the radial projection of aperiodic point sets and study the relation between the resulting
    distribution and properties of the underlying tiling, like symmetry, order and the algebraic type of the
    inflation multiplier.

    \smallskip

    \bigskip
  \end{abstract}
\end{@twocolumnfalse}
]
{
  \renewcommand{\thefootnote}%
    {\fnsymbol{footnote}}
  \footnotetext[1]{email: \texttt{tjakobi@math.uni-bielefeld.de}}
}


\bigskip


\section{Radial Projection Method}\label{sec:radproj}

Given a locally finite point set ${\varLambda \subset \mathbb{R}^{2}}$,
the following procedure is applied (see Fig.~\ref{fig:radproj_exmp} for an example):
\begin{enumerate}[label=(\emph{\alph*})]
  \item\label{radproj:step_vis}{
        For a reference point $x_0 \in \varLambda$ determine the subset $\varLambda'$
        of points visible from this $x_0$. We call a point $p$ \emph{invisible} if there exists
        a $p_0 \in \varLambda$ such that
        \begin{equation}
          \label{eq:visibility_test}
          \exists \; t \in (0,1) \; : \; p_0 = x_0 + t \cdot (p - x_0) \; \text{.}
        \end{equation}}
  \item\label{radproj:step_proj}{
        Select a radius $R > 0$ and project all ${x \in \varLambda' \cap B_R(x_0)}$ to
        $\partial B_R(x_0)$. If $x$ is given as ${(a, b) \in \Real^2}$, then this amounts
        to mapping $x$ to $\arctan(b/a)$.}
  \item\label{radproj:step_sortnorm}{
        The step in \ref{radproj:step_proj} produces a list of distinct (because of the visibility
        condition) angles $\Phi(R)$. Since $\varLambda$ is locally finite, we can sort the list
        \[ \Phi(R) \; = \; \{ \varphi_1, \ldots, \varphi_n \} \]
        in ascending order. We also renormalise the $\varphi_i$ with the factor
        $\frac{n}{2 \pi}$, such that the mean distance between consecutive entries
        becomes $1$.}
  \item\label{radproj:step_meas}{
        Define $d_i := \varphi_{i+1} - \varphi_i$ and consider the discrete
        probability measure
        \begin{equation*}
          \nu_R \; := \; \frac{1}{n-1} \sum_{i=1}^{n-1}{\delta_{d_i}} \; \text{.}
        \end{equation*}}
  \item\label{radproj:step_radlim}{
        Assuming that it exists, the spacing distribution is now obtained by
        taking the limit ${R \rightarrow \infty}$ in the sense of weak convergence
        of measures.}
\end{enumerate}
The choice of the point $x_0$ can be arbitrary and, in general the limit measure
$\nu := \lim_{R \rightarrow \infty}{\nu_R}$ depends on it. For now, we restrict
ourselves to reference points with \emph{high symmetry} (see Fig.~\ref{fig:ab_patch}
and \ref{fig:tt_patch}). Further investigations are needed to decide whether an
averaging over multiple $x_0$ makes more sense here.
\begin{center}
  \includegraphics[width=2.6cm]{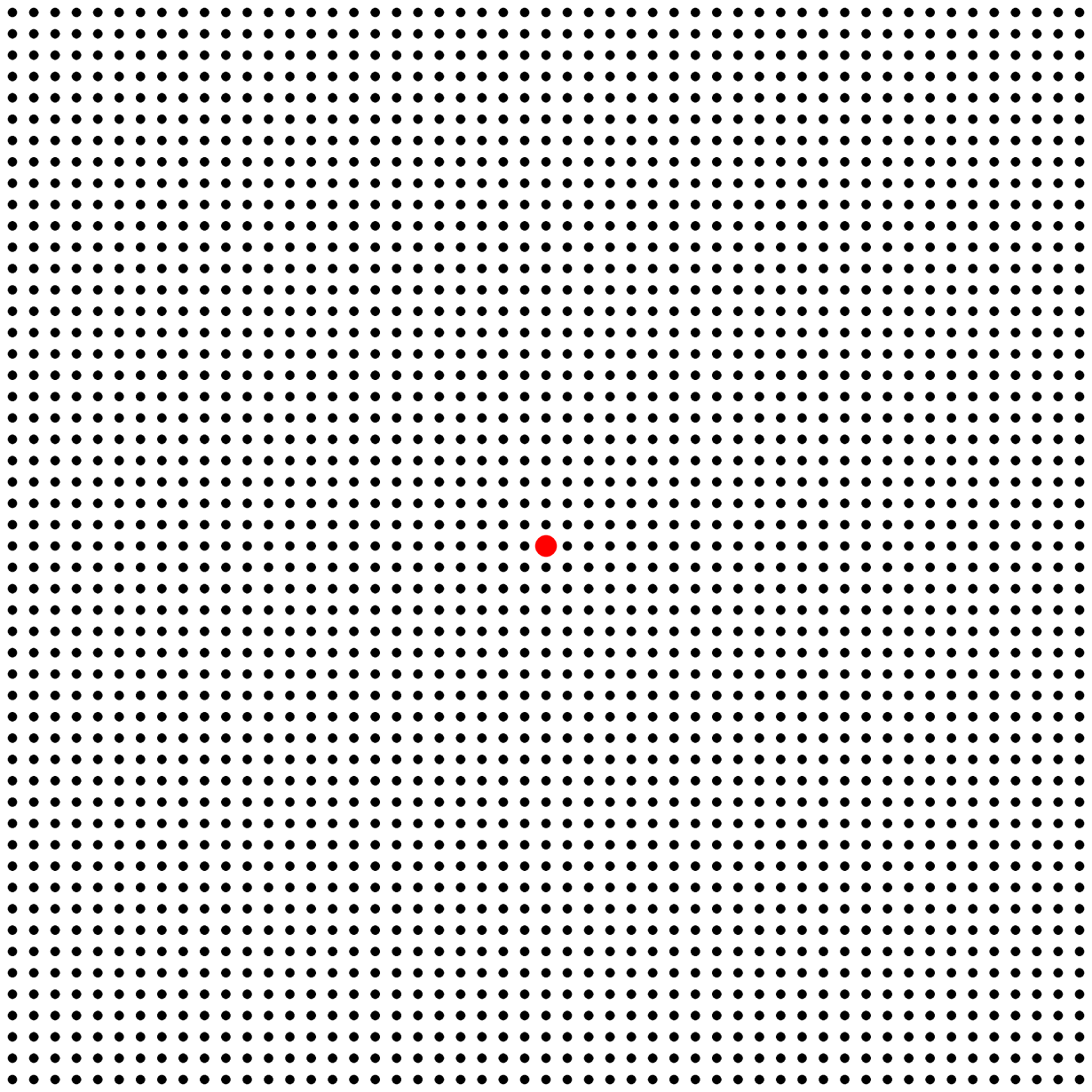}
  \includegraphics[width=2.6cm]{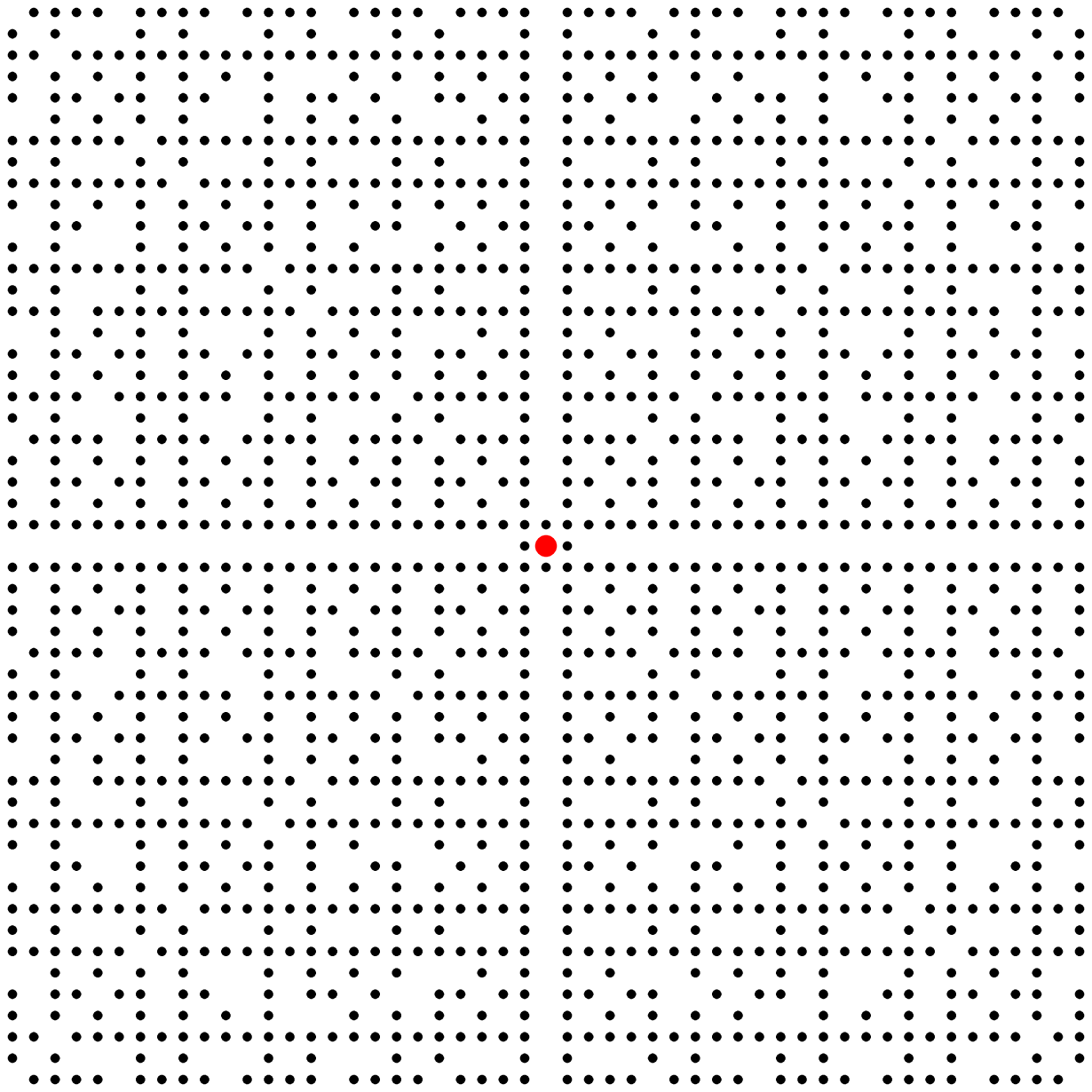}
  \includegraphics[width=2.65cm]{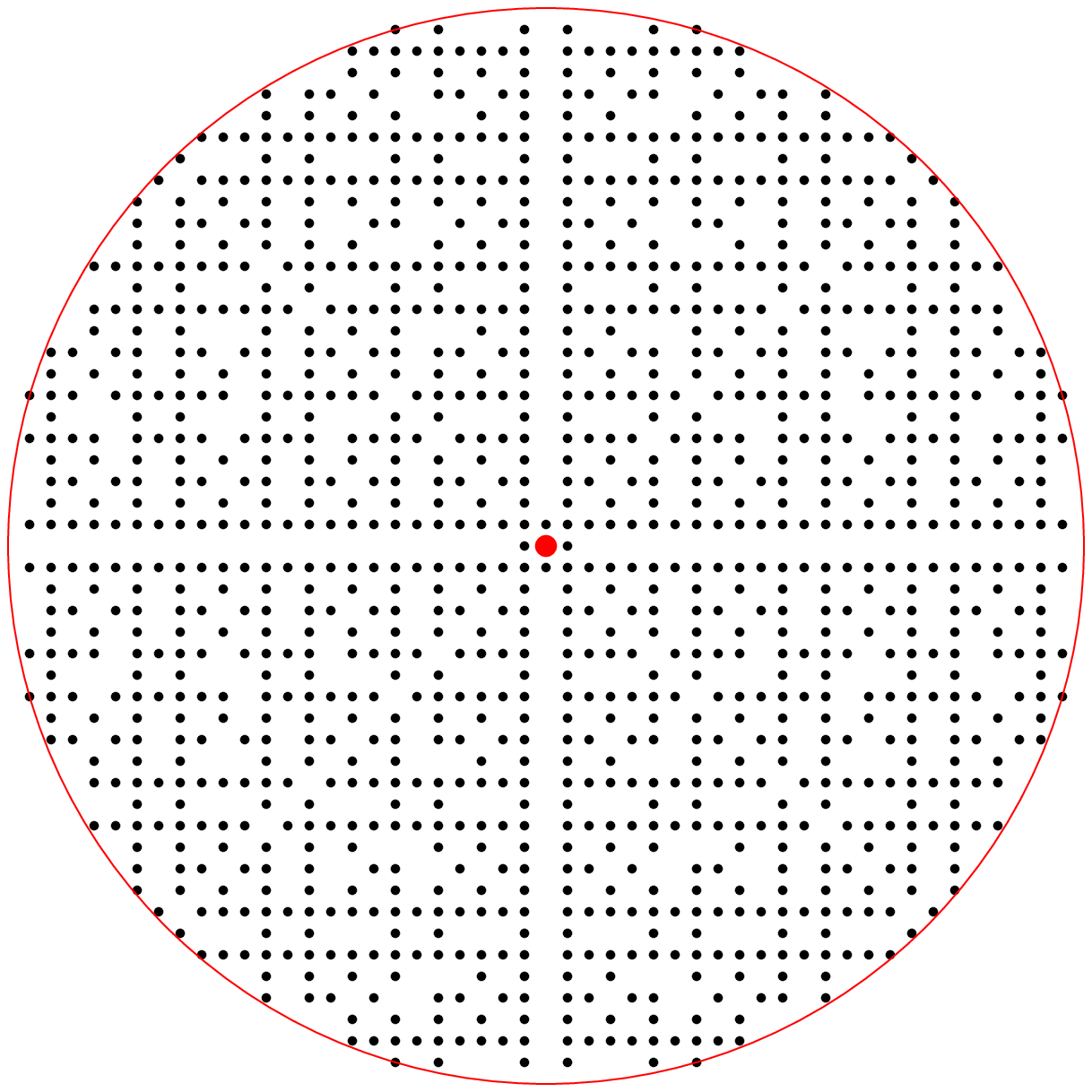}
  \captionof{figure}{Radial projection using the\\
                     example of the $\Integer^2$ lattice.}
  \label{fig:radproj_exmp}
\end{center}
All results (except for the reference cases) are currently strictly numerical and
therefore rely on the computation of large circular patches of the point set and the
determination of visibility. For an introduction to the topic of aperiodic
tilings we refer to \cite{tao13}.

\section{Integer Lattice $\Integer^2$}\label{sec:intlat}

The lattice $\Integer^2$ provides the most ordered case of a planar point set. In this
regard, it represents one reference point set for a potential classification of
order. The set of visible points from the origin is given by
\begin{equation}
  \label{eq:visibility_intlat}
  V_{\Integer^2} \; = \; \{ (a, b) \in \Integer^2 \; : \; \gcd(a, b) = 1 \} \; \text{.}
\end{equation}
In 2000, a closed expression \cite{zaha00} was derived for the limiting measure.
With our setup, the explicit density function reads
\[ g(t) =
  \begin{cases}
    0, & 0 < t < \frac{3}{\pi^2} \text{,} \\
    \frac{6}{\pi^2 t^2} \cdot \log{\frac{\pi^2 t}{3}}, &
      \frac{3}{\pi^2} < t < \frac{12}{\pi^2} \text{,} \\
    \frac{12}{\pi^2 t^2} \cdot \log{\left( 2 \, \Big{/} \! \left( 1 +
      \sqrt{1 - \frac{12}{\pi^2 t}} \; \right) \right)}, &
      t > \frac{12}{\pi^2} \text{,}
  \end{cases}
\]
but in fact the existence also holds for more general expanding regions (which is a
circle here). A Taylor expansion of the tail of $g(t)$ gives
\begin{equation}
  \label{eq:z2_taylor}
  g(1/t) = \frac{36}{\pi^4} t^3 + \frac{162}{\pi^6} t^4 + \mathcal{O}(t^5)
           \; \text{ for } \; t \rightarrow 0_{+} \; \text{,}
\end{equation}
making it obvious that the moments of order ${k \ge 2}$ do not exist. A plot of $g(t)$ is
overlayed over most histograms (see e.g. Fig.~\ref{fig:ab_histo}
and \ref{fig:tt_histo}).\\
For all other cases (apart from the next), only histograms were computed.

\section{Poisson Distributed Points}\label{sec:poisson}

In contrast to the $\Integer^2$ case, the point set generated by a homogeneous spatial
Poisson process gives us the other reference point for our classification. In terms
of order, it represents the case of total disorder. The model is also known as
\emph{complete spatial randomness} (\textsf{CSR}) or \emph{ideal gas} in terms of
physics.\\
Application of the radial projection procedure, with an arbitrary choice of reference
point, yields a $1$-dimensional \textsf{CSR} with intensity ${\lambda = 1}$ (due to the
normalisation) in step \ref{radproj:step_sortnorm}. From there, it is not hard to see
that the density of the radial projection measure is given by
\[ f_{\lambda}(x) \; = \; \begin{cases}
   \lambda \exp(-\lambda x), & x \ge 0,\\
   0,                        & x < 0,
 \end{cases}
\]
since step \ref{radproj:step_meas} just asks the question what the distribution of the
waiting time between jumps of the process is.\\
The density $f_{\lambda}(x)$ (a plot can be seen in Fig.~\ref{fig:gs_histovsexp}) therefore
provides the second explicit function we can test histograms against.

\section{Cyclotomic Model Sets and Visibility}\label{sec:cycloms}

Many aperiodic tilings can not only be realised by inflation of a set of prototiles, but by
projecting a higher-dimensional lattice into $\Real^2$. This is often described as the
\emph{cut and project} method or \emph{model set} description.\\
Our general interest is in \emph{cyclotomic} model sets, since these can be used to construct
$n$-fold symmetric tilings. Also, this description is more suitable for radial projection since
it only gives the tiling vertices and allows precise control over the patch size.\\
The Ammann--Beenker (\textsf{AB}) and the Tübingen triangle (\textsf{TT}) tiling
can both be described as a cyclotomic
\begin{center}
  \includegraphics[width=4.6cm]{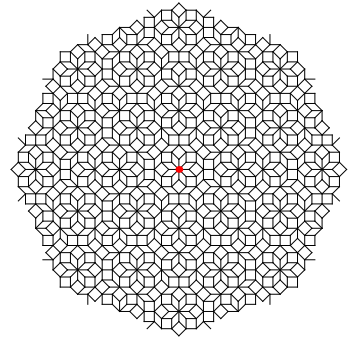}
  \captionof{figure}{\textsf{AB} patch generated via projection.}
  \label{fig:ab_patch}
\end{center}
model set (\textsf{CMS}) and the tiling vertices (the tiling in Fig.~\ref{fig:ab_patch} uses
the triangle/rhombus version) are given by
\begin{align*}
  T_{\textsf{AB}} &= \{ x \in \CyI{8} \; : \; x^{\star} \in W_8 \}, \\
  T_{\textsf{TT}} &= \{ x \in \CyI{5} \; : \; x^{\star} \in W_{10} + \epsilon \},
\end{align*}
where the window $W_8$ is a regular octagon with edge length $1$ centered at the origin.
\begin{center}
  \includegraphics[width=4.6cm]{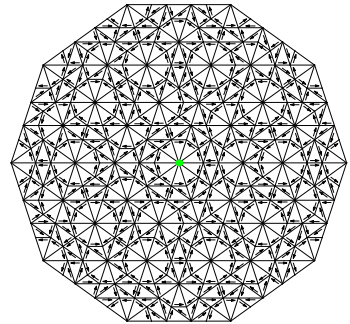}
  \captionof{figure}{\textsf{TT} patch generated via inflation.}
  \label{fig:tt_patch}
\end{center}
For the \textsf{TT} tiling the window $W_{10}$ is a decagon, here with edge length $\sqrt{(\tau + 2)/5}$
($\tau$ denoting the golden mean). For the orientations of the windows, see Fig.~\ref{fig:ab_phyvsint}.
The map $\star$ is given by the extension of $\zeta_8 \mapsto \zeta_8^3$ in the \textsf{AB} case, and
$\zeta_5 \mapsto \zeta_5^2$ in the \textsf{TT} one (${\zeta_n = \exp(2 \pi i / n)}$).
\begin{center}
  \includegraphics[width=5.8cm]{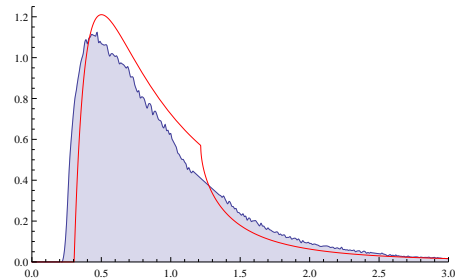}
  \captionof{figure}{Spacing distribution of a large \textsf{AB} patch.}
  \label{fig:ab_histo}
\end{center}
Another aspect of choosing this description is that the visibility test in Eq.\eqref{eq:visibility_test}
is computationally expensive compared to the local test \eqref{eq:visibility_intlat}
of the $\Integer^2$ case. It turns out that a \textsf{CMS} also admits a similar description.
For example, the visible points of the \textsf{AB} tiling with the reference point at
the origin are
\[ V_{\textsf{AB}} = \{ x \in T_{\textsf{AB}} \; : \; \lambda_{\text{sm}} x^{\star}
   \notin W_8 \text{ and } x \text{ is coprime} \}, \]
with ${\lambda_{\text{sm}} = 1 + \sqrt{2}}$, the silver mean. The $\Integer$-module in this
case can be decomposed into $\CyI{8} = \Integer[\sqrt{2}] \oplus \Integer[\sqrt{2}] \cdot \zeta_8$
and coprime in this context means that for $\CyI{8} \ni x = x_1 + x_2 \cdot \zeta_8$ the
$\gcd(x_1, x_2)$ is a unit in $\Integer[\sqrt{2}]$.
\begin{center}
  \includegraphics[width=5.8cm]{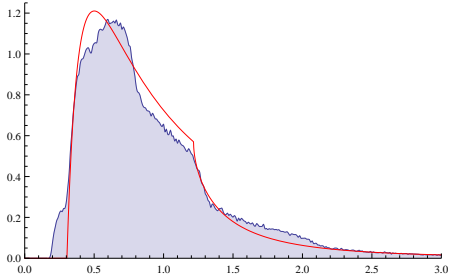}
  \captionof{figure}{Spacing distribution of a large \textsf{TT} patch.}
  \label{fig:tt_histo}
\end{center}
The same description can be given for the \textsf{TT} tiling and reads
\[ V_{\textsf{TT}} \; = \; \{ x \in T_{\textsf{TT}} \; : \; \tau x^{\star} \notin
   W_{10} + \epsilon \text{ and } x \text{ is coprime} \} \; \text{,} \]
with $\epsilon$ a small shift. Coprimality is defined as in the \textsf{AB}
case, except that the module is $\Integer[\tau]$ this time. The simplicity of
these local visibility tests also depends on the order of
\begin{center}
  \includegraphics[width=5.0cm]{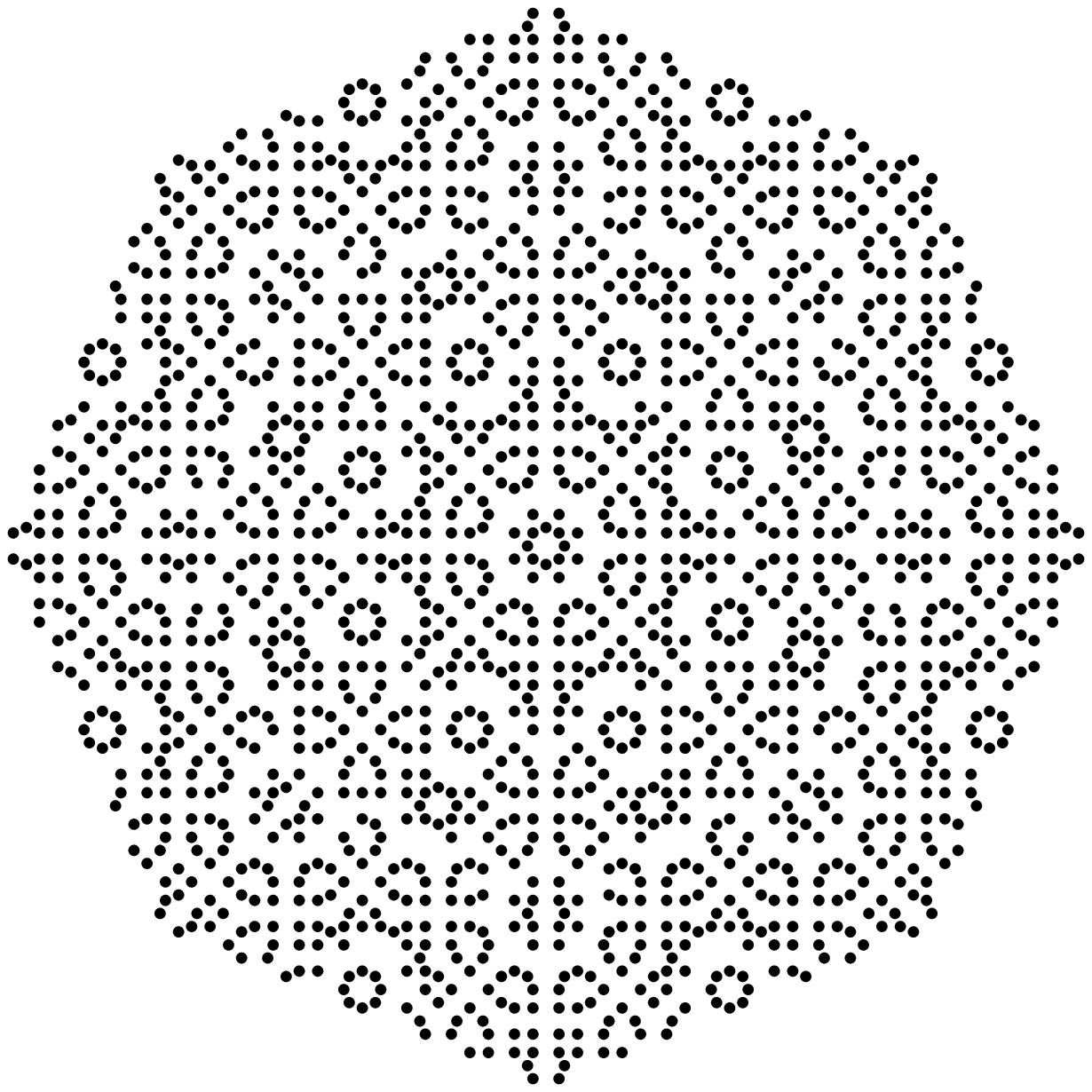} \hspace{4.1ex}
  \raisebox{1.15cm}{\includegraphics[width=2.6cm]{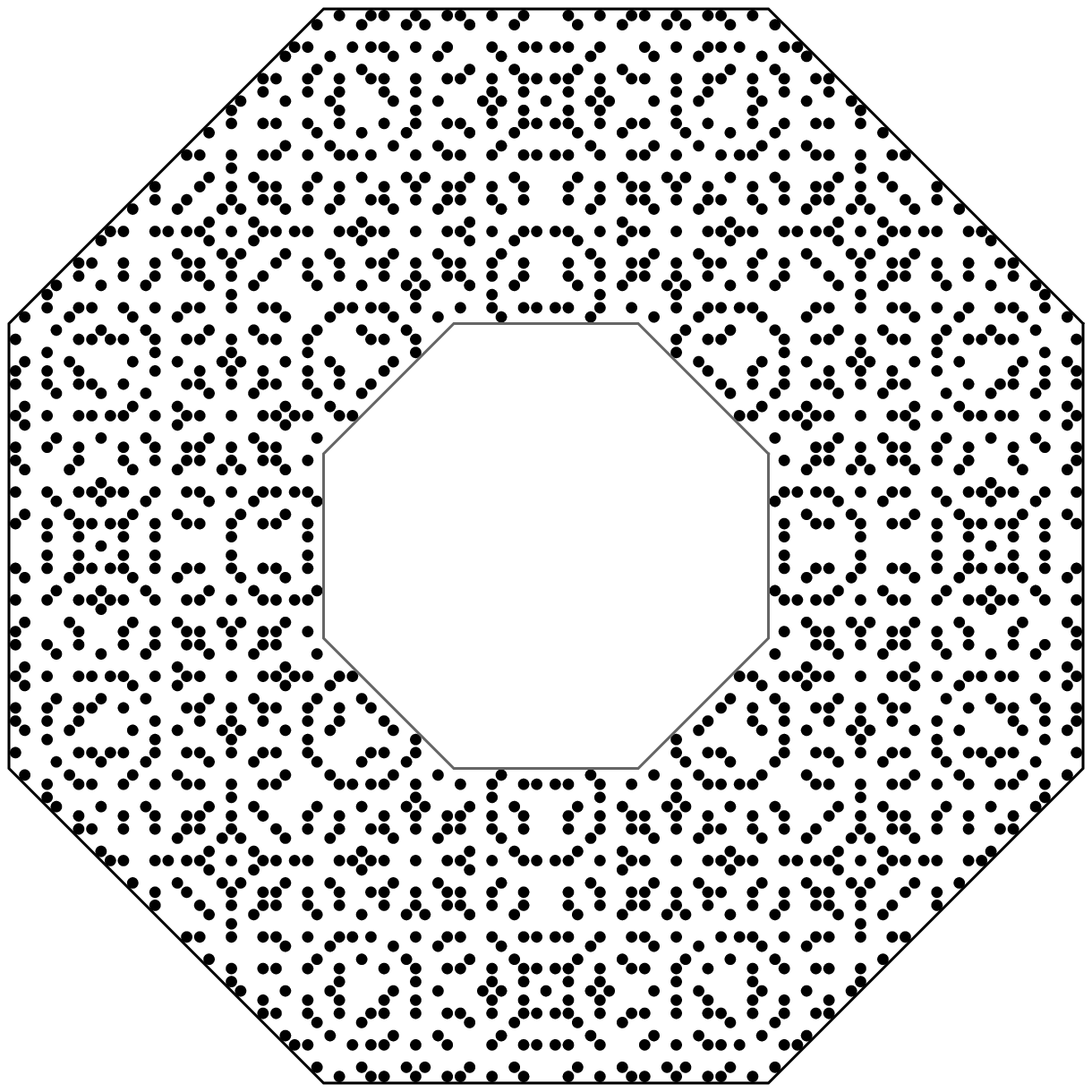}}
  \captionof{figure}{Visible vertices of the $8$-fold \textsf{AB} tiling\\
                     (left: direct space, right: internal space).}
  \label{fig:ab_phyvsint}
\end{center}
the underlying cyclotomic field, here a prime power.

\section{A $12$-fold Cyclotomic Case}\label{sec:cyclotwelve}

Another tiling that is given by a CMS is the \emph{Gähler shield} (\textsf{GS}) tiling. The
vertices are
\[ T_{\textsf{GS}} \; = \; \{ x \in \CyI{12} \; : \; x^{\star} \in W_{12} + \epsilon \} \; \text{,} \]
with a dodecagon $W_{12}$ (edge length $1$) and the map $\star$ given by the extension of
$\zeta_{12} \mapsto \zeta_{12}^5$. The order of the cyclotomic field leads to a slightly more
involved local visibility test
\begin{align*}
  V_{\textsf{GS}} \; = \; &\{ x \in T_{\textsf{GS}} \; : \; n(x) = 1 \wedge
                            \lambda_1 x^{\star} \notin W_{12} + \epsilon \} \cup {} \\
                          &\{ x \in T_{\textsf{GS}} \; : \; n(x) = 2 \wedge
                            \lambda_2 x^{\star} \notin W_{12} - \epsilon \} \; \text{.}
\end{align*}
The function $n$ decomposes a $x \in \CyI{12}$ into the direct-sum representation
\[ \Integer[\sqrt{3}] \oplus \Integer[\sqrt{3}] \cdot \zeta_{12} \; \text{ with } \;
   \lambda_{12} := 2 + \sqrt{3} \]
and then computes $\abs{N(\gcd(x_1, x_2))}$, the absolute value of the algebraic norm of the $\gcd$ of the components.
The two \emph{rescaling factors} for the visibility test are given by
\[ \lambda_1 := \sqrt{\lambda_{12} \cdot 2} \quad \text{ and } \quad
   \lambda_2 := \sqrt{\lambda_{12} / 2} \; \text{.} \]
The first part of the set $V_{\textsf{GS}}$ (indicated in \textcolor{lgrey}{\textbf{grey}} in
Fig.~\ref{fig:gs_phyvsint}) comprises the already known coprime elements.
\begin{center}
  \includegraphics[width=5.0cm]{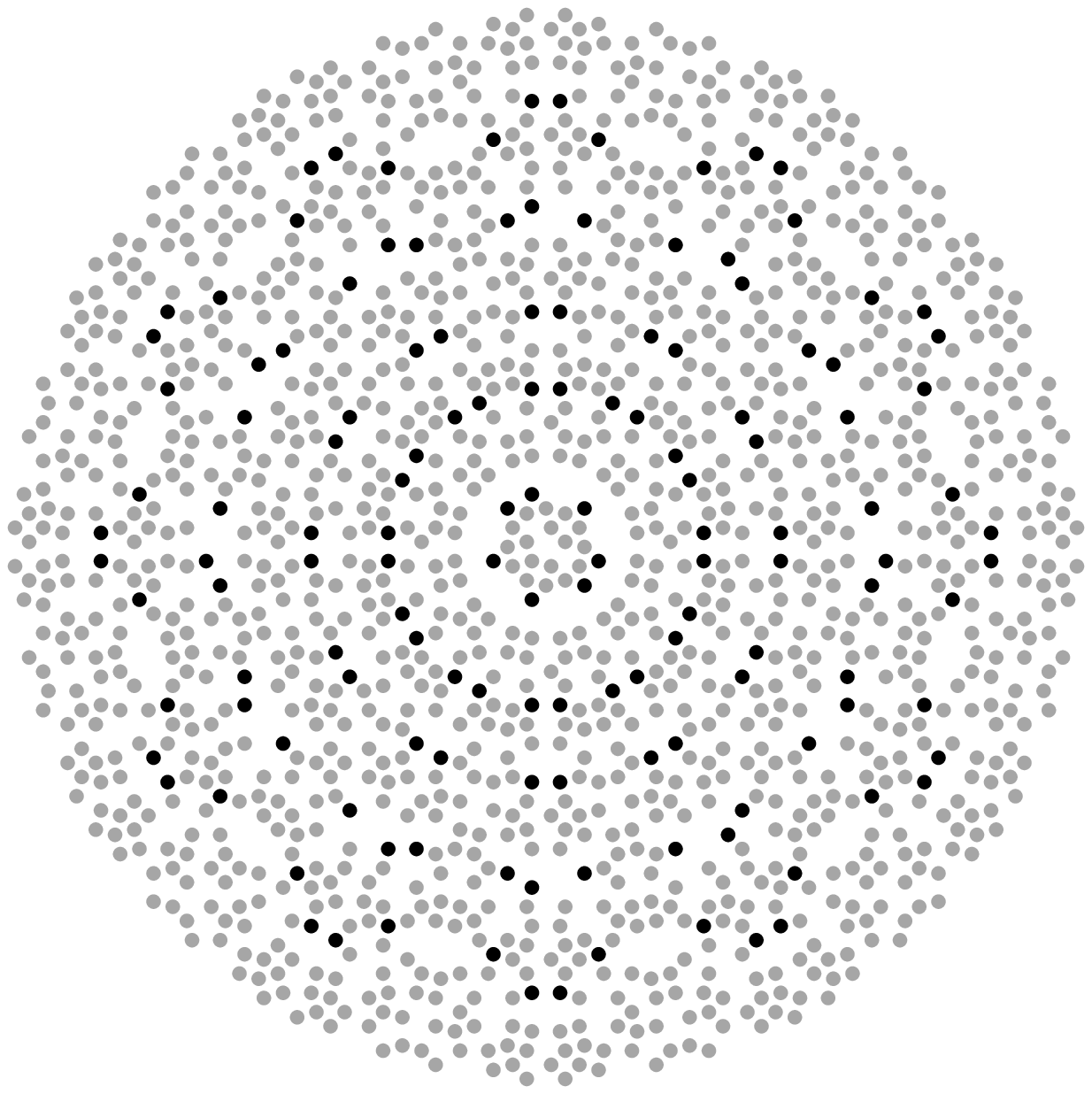} \hspace{4.1ex}
  \raisebox{1.05cm}{\includegraphics[width=2.6cm]{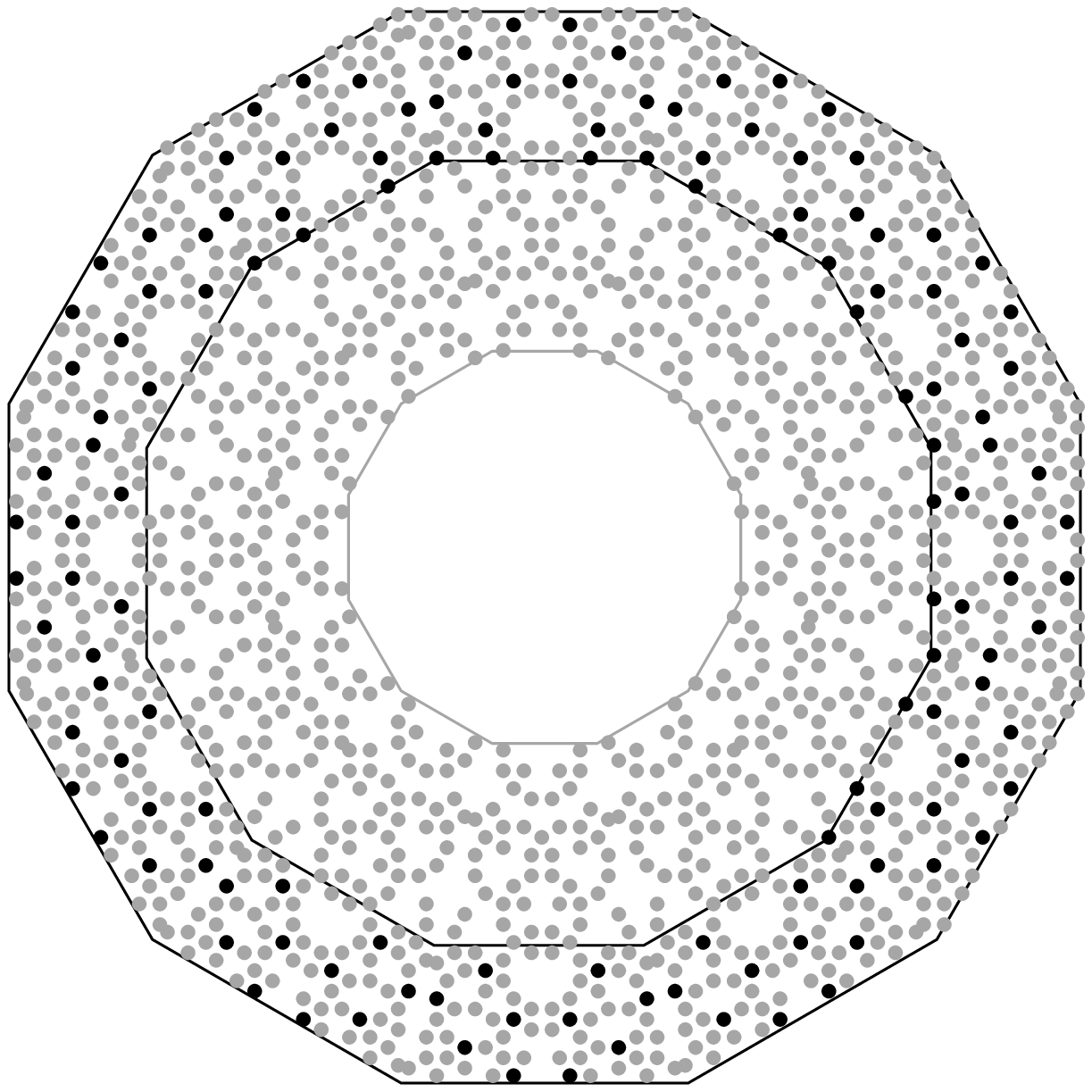}}
  \captionof{figure}{Visible points of a $12$-fold \textsf{GS} tiling\\
                     (left: direct space, right: internal space)}
  \label{fig:gs_phyvsint}
\end{center}
The \textbf{second part} is exceptional, and its existence is linked to the order
$12$ being a composite number.
\begin{center}
  \includegraphics[width=5.9cm]{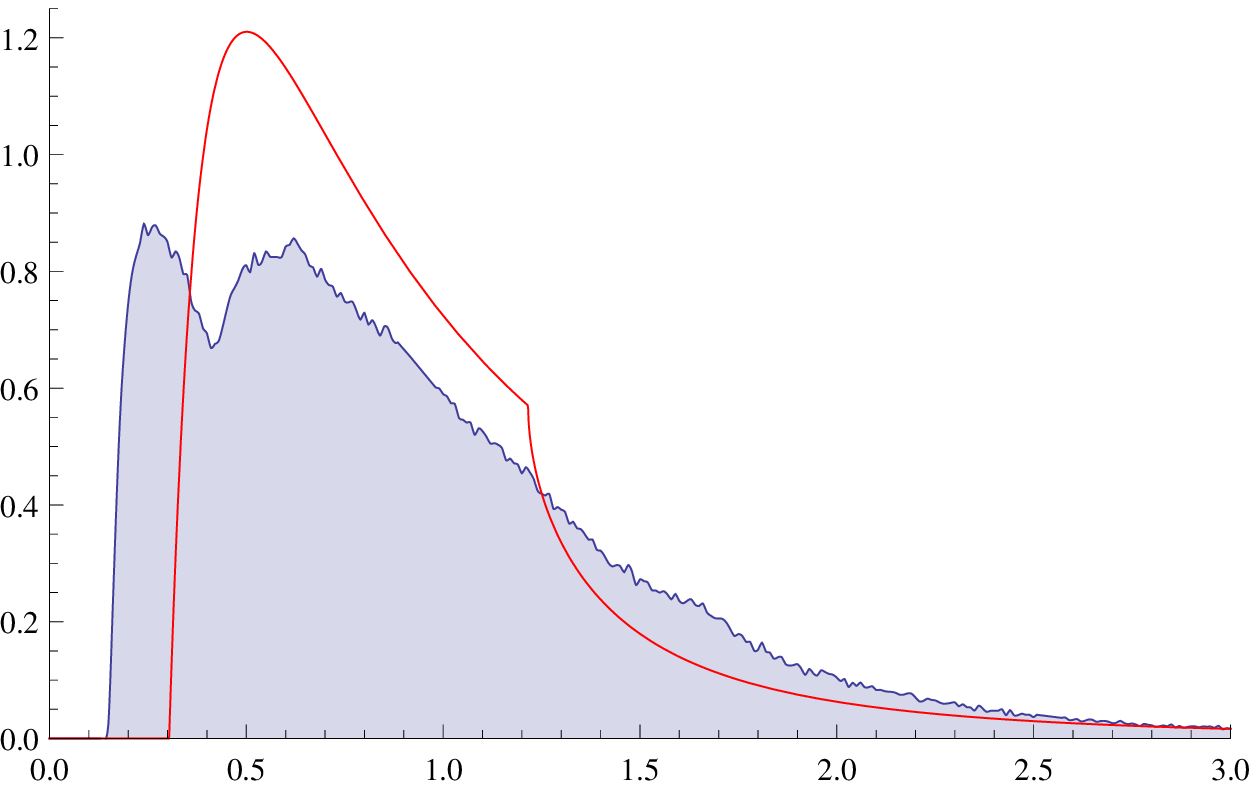}
  \captionof{figure}{Spacing distribution of a large \textsf{GS} patch.}
  \label{fig:gs_histovsz2}
\end{center}
Even though Fig.~\ref{fig:gs_histovsz2} deviates a lot more from the $\Integer^2$ distribution
(compared to Fig.~\ref{fig:ab_histo} and \ref{fig:tt_histo}), all considered cyclotomic cases
exhibit the special threefold structure (gap, bulk and tail). Let us consider a different tiling
\begin{center}
  \includegraphics[width=5.9cm]{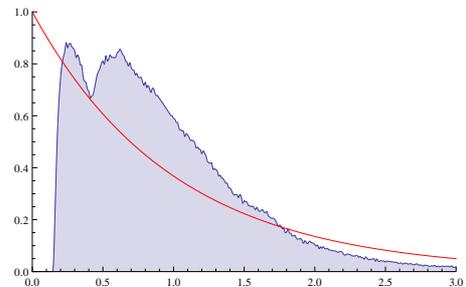}
  \captionof{figure}{\textsf{GS} spacing distribution\\compared against Poisson.}
  \label{fig:gs_histovsexp}
\end{center}
with stronger spatial fluctuations.

\section{A Non-Pisot Case / Lan\c{c}on--Billard}\label{sec:nonpisot}

The Lan\c{c}on-Billard (\textsf{LB}) tiling is an example of a \emph{chiral} (the tiles only
appear in one chirality) inflation tiling with
\begin{center}
  \includegraphics[width=2.55cm]{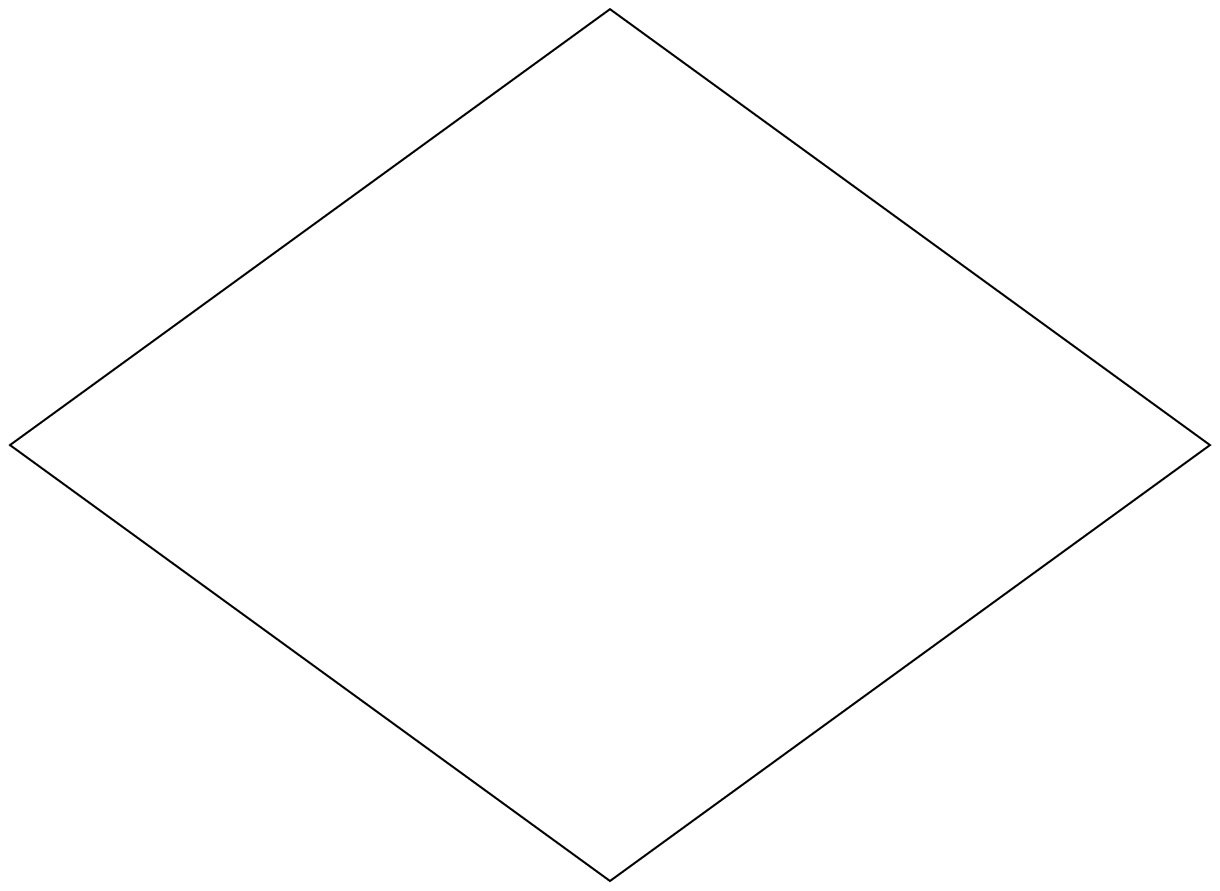} \hspace{0.5ex}
  \raisebox{0.84cm}{$\xrightarrow{\hspace{4ex}}$} \hspace{0.5ex}
  \includegraphics[width=2.55cm]{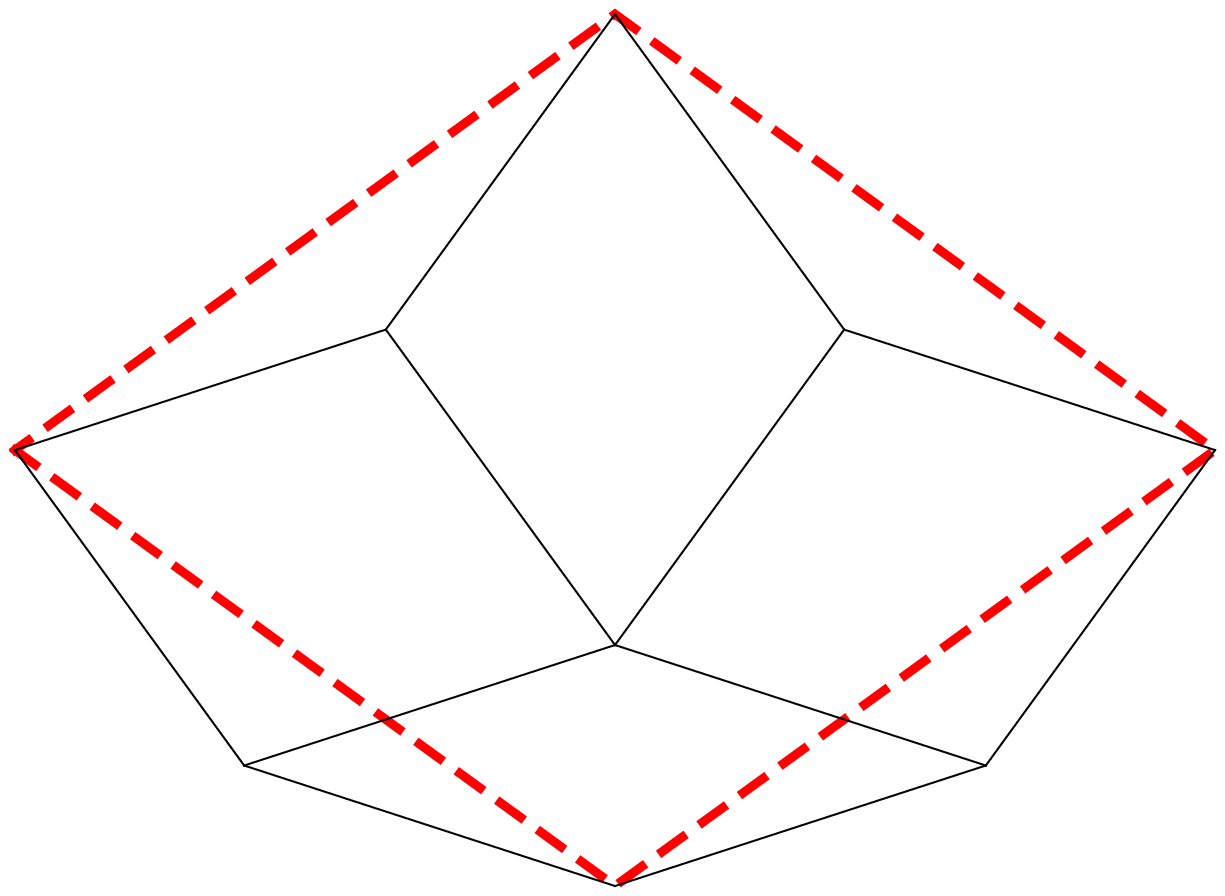}
  \captionof{figure}{Tile A maps to 3$\times$A and 1$\times$B.}
  \label{fig:ch_infl_a}
\end{center}
a multiplier $\lambda_{\textsf{LB}} = \sqrt{(5 + \sqrt{5})/2}$, which
is a non-Pisot number. The inflation rules (a \emph{stone} inflation is only
\begin{center}
  \includegraphics[width=2.55cm]{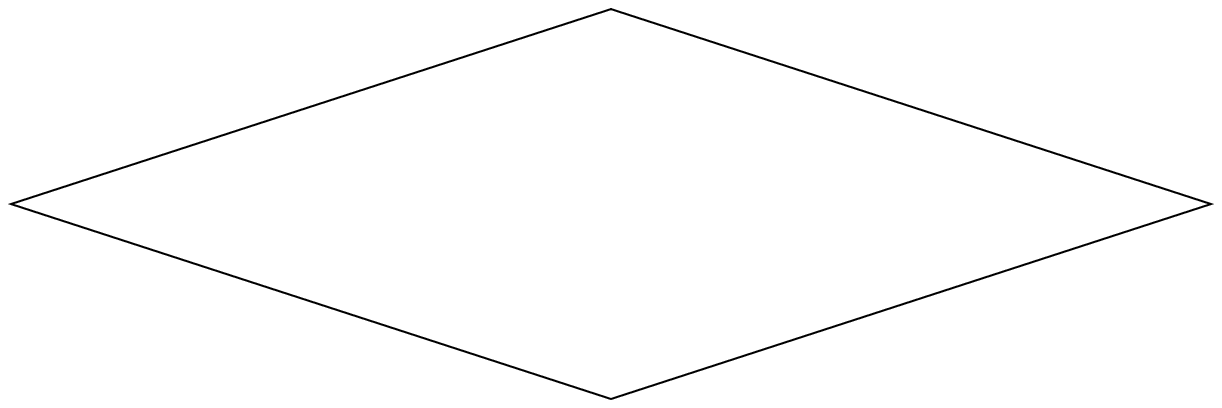} \hspace{0.5ex}
    \raisebox{0.33cm}{$\xrightarrow{\hspace{4ex}}$} \hspace{0.5ex}
    \includegraphics[width=2.55cm]{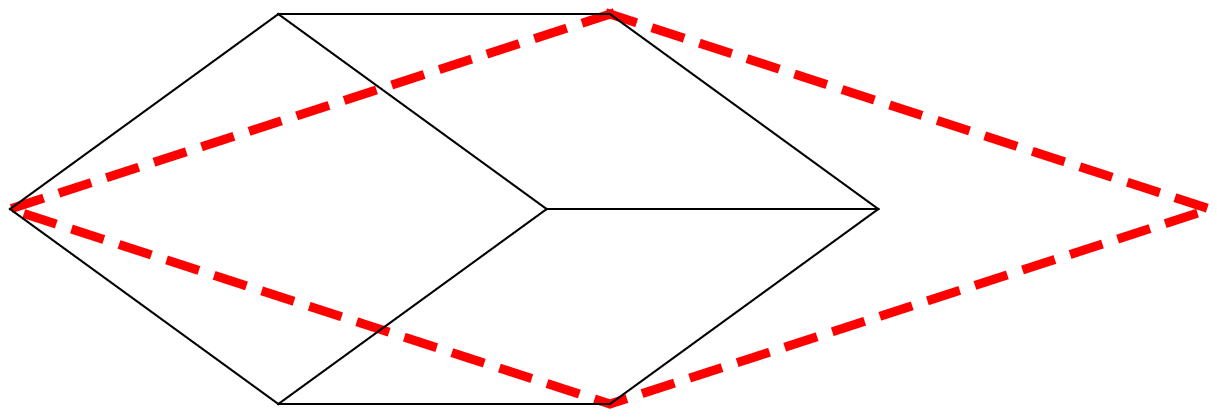}\vspace*{0.4ex}
    \captionof{figure}{Tile B maps to 1$\times$A and 2$\times$B.}
    \label{fig:ch_infl_b}
\end{center}
possible with a fractal tile boundary) generate a very irregular tiling.
\begin{center}
  \includegraphics[width=4.8cm]{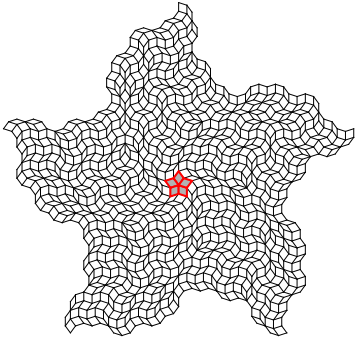}
  \captionof{figure}{\textsf{LB} $5$-fold symmetric patch.}
  \label{fig:ch_circ4}
\end{center}
The resulting distribution shows that the radial projection is sensitive
to this irregularity.
\begin{center}
  \includegraphics[width=5.8cm]{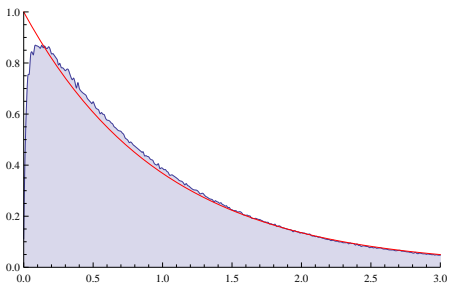}
  \captionof{figure}{Spacing distribution of a large \textsf{LB} patch.}
  \label{fig:ch_histo}
\end{center}
In fact, on this level, the \textsf{LB} tiling is almost indistinguishable
from the Poisson case.

\section{Some Additional Examples}\label{sec:morexmp}

We have already seen in Secs.~\ref{sec:cycloms} and \ref{sec:cyclotwelve} that the
radial projection reacts to the order of symmetry of the tiling. The vertices of
the \emph{chair} tiling in Fig.~\ref{fig:histo_other} are a subset of
\begin{center}
  \includegraphics[width=4.2cm]{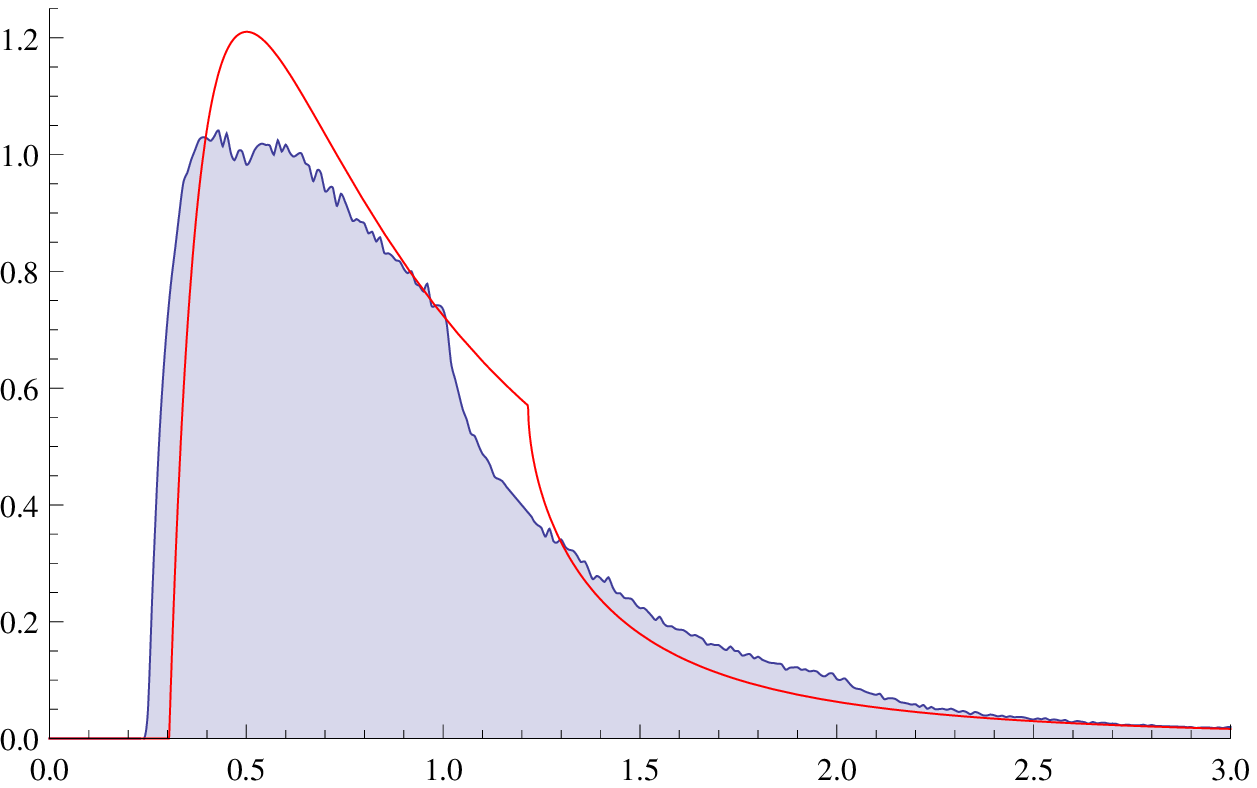} \hspace{0.8ex}%
  \includegraphics[width=4.2cm]{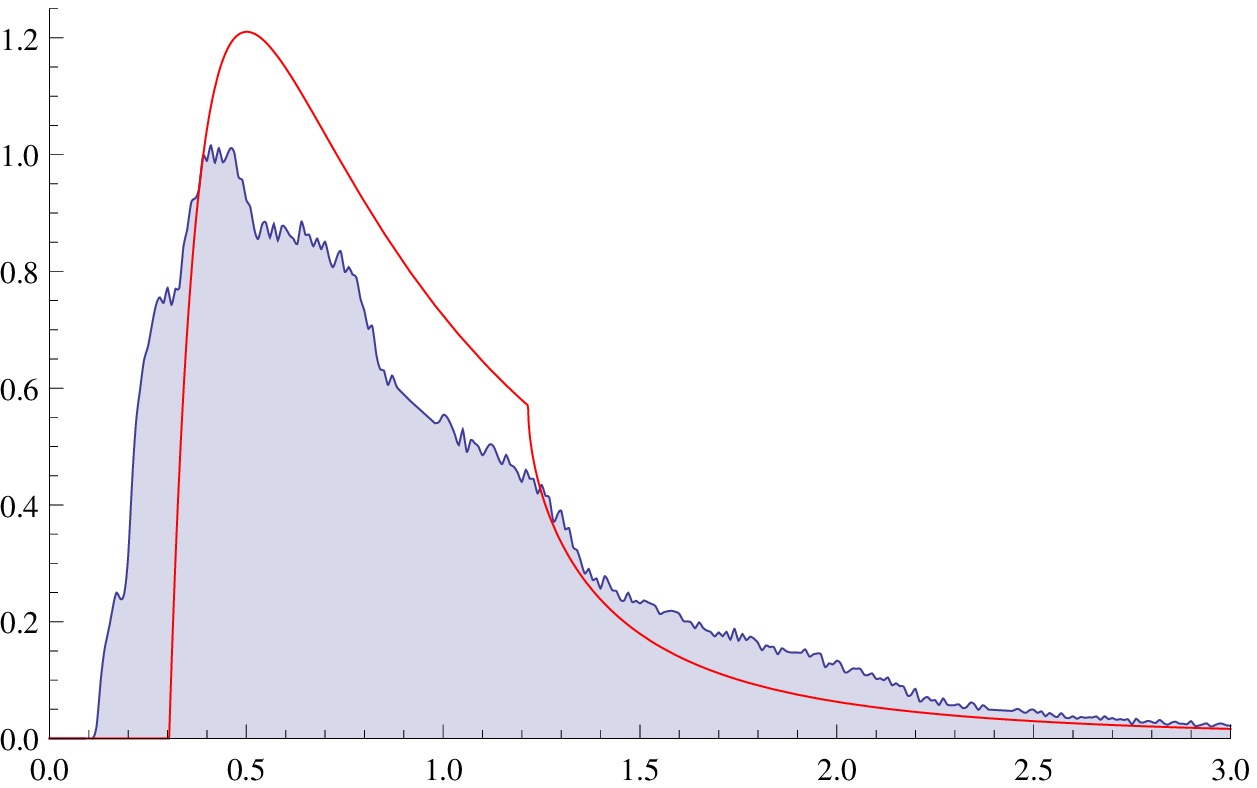}
  \captionof{figure}{Spacing distribution of the \emph{chair} (left)\\
                     and \emph{rhombic Penrose} (right) tiling.}
  \label{fig:histo_other}
\end{center}
$\Integer^2$, but a different type of visibility condition
has to be applied. The bulk section notices these changes, and
\begin{center}
  \includegraphics[width=4.2cm]{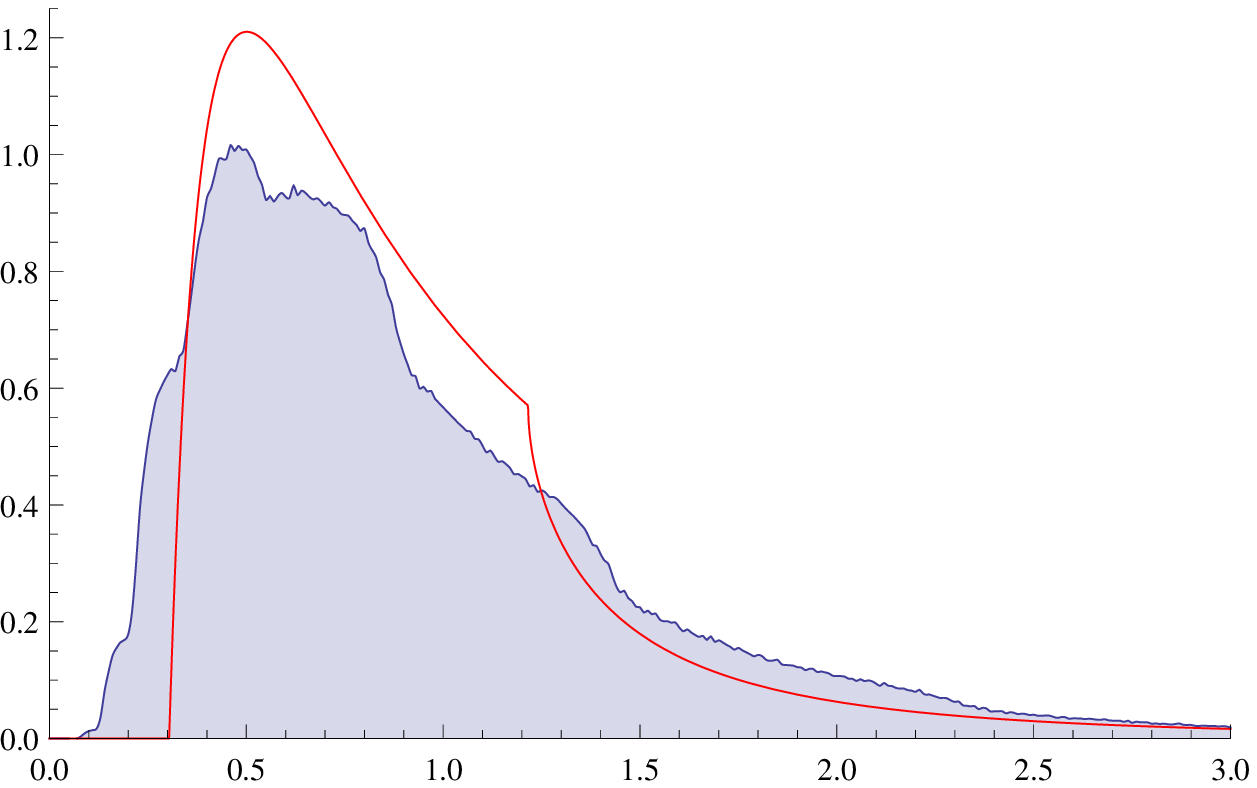} \hspace{0.8ex}%
  \includegraphics[width=4.15cm]{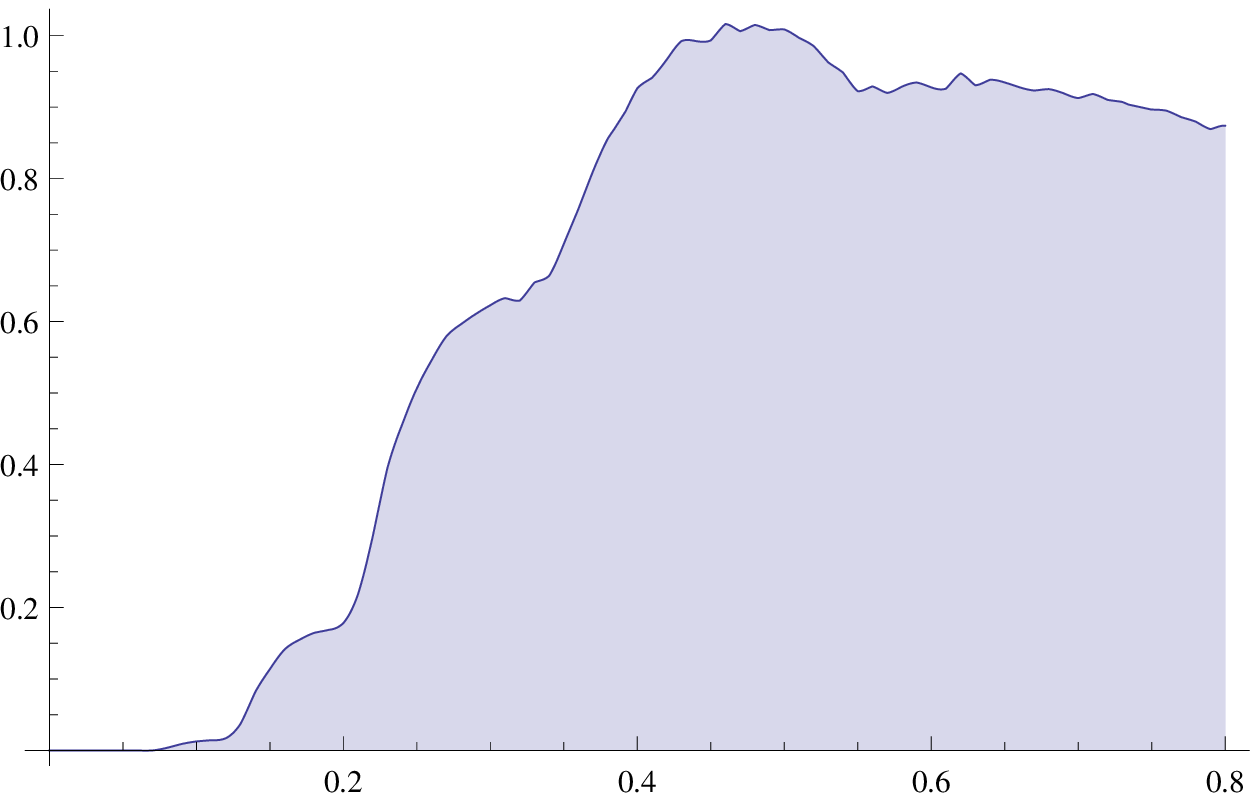}
  \captionof{figure}{Spacing distribution of the \emph{Penrose-Robinson} tiling
                     and zoom into the bulk.}
  \label{fig:prt_and_zoom}
\end{center}
shows a lot of structure in general.

\section{Tail Decay Behaviour}\label{sec:taildecay}

The expansion in Eq. \eqref{eq:z2_taylor} gives us an idea how the tail decay of the
distribution behaves. A power law fit ($c_k$ the coefficient of $t^k$) applied to
the numerical data
\begin{center}
  \begin{tabular}{c|c|c|c|c}\toprule
    tiling       & gap size  & $c_3$  & $c_4$  & $e$  \\
    \midrule
    $\Integer^2$ & 0.304     & 0.369  & 0.168  & ---  \\
    \textsf{AB}  & 0.222     & 0.248  & 0.496  & 2.79 \\
    \textsf{TT}  & 0.182     & 0.239  & 0.513  & 2.60 \\
    \textsf{GS}  & 0.152     & 0.232  & 0.547  & 4.75 \\
    \bottomrule
  \end{tabular}
  \captionof{table}{Statistical data generated from the\\
                    radial projection ($e = $ error $\cdot 10^{-10}$).}
  \label{tab:histo_stats}
\end{center}
indicates that at least the \textsf{CMS} cases also show a similar kind of decay.

\section*{Acknowledgement}

The author wishes to thank Michael Baake, Markus Moll and Christian Huck for
helpful discussions. This work is supported by the German Research Foundation
(DFG) via the Collaborative Research Centre (CRC 701) at the faculty of
Mathematics, University of Bielefeld.

\end{document}